\documentclass{article}
\usepackage{amsfonts}
\usepackage{amsmath}
\usepackage{amsmath,amssymb,latexsym,color}
\usepackage[mathscr]{eucal}
\newtheorem{thm}{Theorem}[section]

\newtheorem{lemma}[thm]{Lemma}
\newtheorem{cor}[thm]{Corollary}

\newtheorem{example}{Example}[section]
\newtheorem{defin}[thm]{Definition}

\newtheorem{remark}{Remark}[section]

\newcommand{\qed}{\hfill\Box\medskip}
\usepackage{CJK}
\begin{document}
\begin{CJK*}{GBK}{song}

\renewcommand{\baselinestretch}{1.3}
\title{\bf The Terwilliger algebra of the incidence graphs of Johnson geometry}

\author{
Qian Kong\quad Benjian Lv\quad
Kaishun Wang\footnote{Corresponding author. E-mail address: wangks@bnu.edu.cn} \\
{\footnotesize \em  Sch. Math. Sci. {\rm \&} Lab. Math. Com. Sys.,
Beijing Normal University, Beijing, 100875,  China} }
\date{}
\maketitle

\begin{abstract}
Levstein and Maldonado [F. Levstein, C. Maldonado, The Terwilliger
algebra of the Johnson schemes, Discrete Mathematics 307 (2007)
1621--1635] computed the Terwilliger algebra of the Johnson scheme
$J(n,m)$ when $3m\leq n$. In this paper, we determine the
Terwilliger algebra of the incidence graph $J(n,m,m+1)$ of Johnson
geometry when $3m\leq n$, give two bases of this algebra, and
calculate its dimension.



\end{abstract}

\section{Introduction}

Let $\Gamma=(X,R)$ denote a simple connected graph with the vertex
set $X$ and the edge set $R$. For vertices $x$ and $y$,
$\partial(x,y)$ denotes the \emph{distance} between $x$ and $y$,
i.e., the length of a shortest path connecting $x$ and $y$. Fix a
vertex $x\in X$. Let $D(x):=\max\{\partial(x,y)\mid y\in X\}$ denote
the \emph{diameter with respect to} $x$. For each $i\in
\{0,1,\ldots,D(x)\}$, let $\Gamma_i(x)=\{y\in X\mid
\partial(x,y)=i\}$ and define $E_i^*=E_i^*(x)$ to be the diagonal matrix in
Mat$_X(\mathbb{C})$ with $yy$-entry
$$
(E_i^*)_{yy}=\left\{ \begin{array}{ll}
1,& \textrm{if $y\in \Gamma_i(x)$},\\
0,& \textrm{otherwise}.\\
\end{array} \right.
$$
Let $\mathcal {T}=\mathcal {T}(x)$ be the subalgebra of
Mat$_X(\mathbb{C})$ generated  by the adjacency matrix $A$ of
$\Gamma$ and $E_0^*,E_1^*,\ldots,E_{D(x)}^*$. Then $\mathcal {T}$ is
called the \emph{Terwilliger algebra} of $\Gamma$ with respect to
$x$. Let $V=\mathbb{C}^X$ denote the vector space over the complex
number field $\mathbb{C}$ consisting of column vectors whose
coordinates are indexed by $X$. A \emph{$\mathcal {T}$-module} is
any subspace $W\subseteq V$ such that $\mathcal {T}W\subseteq W$. We
call a nonzero $\mathcal {T}$-module $W$ \emph{irreducible} if it
does not properly contain a nonzero $\mathcal {T}$-module. An
irreducible $\mathcal {T}$-module $W$ is \emph{thin} if $\dim
E_i^*W\leq 1$ for every $i$, and the graph $\Gamma$ is said to be
\emph{thin with respect to $x$} if every irreducible $\mathcal
{T}(x)$-module is thin.

Terwilliger \cite{ter1, ter2, ter3} initiated the study of the
Terwilliger algebra of association schemes, which has been used to
study group schemes \cite{s5a5, bannai1}, strongly regular graphs
\cite{sr}, bipartite and almost bipartite $P$- and $Q$-polynomial
association schemes \cite{bpq, abpq}, $2$-homogeneous bipartite
distance-regular graphs \cite{curtin}, the Hypercube \cite{go}, the
Hamming schemes \cite{hdq} and the Johnson schemes \cite{johnson},
etc.

Let $\Omega$ be a set of cardinality $n$ and let ${\Omega \choose
i}$ denote the set of all $i$-subsets of $\Omega$. The
\emph{incidence graph $J(n,m,m+1)$ of the Johnson geometry} is a
bipartite graph with a bipartition ${\Omega \choose m}\cup {\Omega
\choose m+1}$, where $y\in {\Omega \choose m}$ and $z\in {\Omega
\choose m+1}$ are adjacent if and only if $y\subseteq z$. It is
known that $J(n,m,m+1)$ is \emph{distance-biregular}(see \cite
{bcn}).

Levstein and Maldonado \cite{johnson} determined the Terwilliger
algebra of the Johnson scheme $J(n,m)$ when $3m\leq n$. Motivated by
this result, in this paper we shall determine the Terwilliger
algebra of $J(n,m,m+1)$ with respect to $x\in {\Omega \choose m}$
when $n\geq 3m$.

This paper is organized as follows. In Section 2, we introduce the
intersection matrices and give some useful identities. In Section 3,
we determine the Terwilliger algebra of $J(n,m,m+1)$, and show
$J(n,m,m+1)$ is thin with respect to $x$. In Section 4, we give two
bases of the Terwilliger algebra and compute its dimension.

\section{Intersection matrices}

In this section we first introduce the inclusion matrices of a set,
then discover the relationship between the adjacency matrix of
$J(n,m,m+1)$ and the inclusion matrices, and give some identities
for intersection matrices.

The following lemma is useful.

\begin{lemma}\label{1}
Let $J(n,m,m+1)$ be the incidence graph of Johnson geometry with a
bipartition ${\Omega \choose m}\cup {\Omega \choose m+1}$. Pick
$x\in {\Omega \choose m}$. Then $\partial(x,z)=2i$ if and only if
$|z|=m$ and $|x\cap z|=m-i$; $\partial(x,z)=2i+1$ if and only if
$|z|=m+1$ and $|x\cap z|=m-i$. Furthermore, when $n\geq 2m+1$ we
have $D(x)=2m+1$.
\end{lemma}
\textbf{Proof.} Immediate from \cite [Lemma 2.2 (1)(3)] {hiraki}.
$\qed$

Fix $x\in {\Omega \choose m}$. We then consider the adjacency matrix
$A$ of $J(n,m,m+1)$ as a block-matrix with respect to the partition
$\{x\}\cup\Gamma_1(x)\cup\cdots\cup\Gamma_{2m+1}(x)$. In order to
describe the blocks of $A$, we need to introduce the inclusion
matrices.

Let $V$ be a set of cardinality $v$. The \emph{inclusion matrix}
$W_{i,j}(v)$ is a $(0,1)$-matrix whose rows and columns are indexed
by ${V \choose i}$ and ${V \choose j}$, respectively, with the
$yz$-entry defined by
$$
(W_{i,j}(v))_{yz}=\left\{ \begin{array}{ll}
1,& \textrm{if $y\subseteq z$},\\
0,& \textrm{otherwise}.\\
\end{array} \right.
$$
Observe that
\begin{eqnarray}\label{3}
W_{i,j}(v)W_{j,k}(v)={k-i \choose j-i}W_{i,k}(v).
\end{eqnarray}

Let $A_{i,j}$ be the submatrix of $A$ with rows indexed by
$\Gamma_i(x)$ and columns indexed by $\Gamma_j(x)$.

\begin{lemma}\label{4}
Let $I_{{v \choose k}}$ be the identity matrix of size ${v \choose
k}$. Then
\begin{eqnarray}
\label{5}{}&&A_{i,j}=0 \quad (0\leq i\leq j \leq 2m+1\ \
\textrm{and} \ \ i\neq j-1),{}\\
\label{6}{}&&A_{2i,2i+1}=I_{{m \choose m-i}}\otimes W_{i,i+1}(n-m) \quad (0\leq i\leq m),\\
\label{7}{}&&A_{2i+1,2i+2}=(W_{m-i-1,m-i}(m))^\emph{t}\otimes
I_{{n-m \choose i+1}} \quad (0\leq i\leq m-1),
\end{eqnarray}
where ``$\otimes$'' denotes the Kronecker product of matrices.
\end{lemma}
\textbf{Proof.} (\ref{5}) is directed.

Pick $y\in\Gamma_{2i}(x)$, $z\in\Gamma_{2i+1}(x)$. By Lemma \ref{1}
we have $|y|=m$, $|z|=m+1$, $|x\cap y|=|x\cap z|=m-i$. Suppose
$y=\alpha_{m-i}\beta_i:=\alpha_{m-i}\cup\beta_i$,
$z=\alpha_{m-i}'\beta_{i+1}'$, where $\alpha_{m-i}$ and
$\alpha_{m-i}'\in{x \choose m-i}$, while $\beta_i\in{\Omega\setminus
x \choose i}$ and $\beta_{i+1}'\in{\Omega\setminus x \choose i+1}$.
Then
$$
(A_{2i,2i+1})_{yz}=(I_{m \choose m-i}\otimes
W_{i,i+1}(n-m))_{yz}=\left\{ \begin{array}{ll}
1,& \textrm{if $\alpha_{m-i}=\alpha_{m-i}'$ and $\beta_i\subseteq \beta_{i+1}'$},\\
0,& \textrm{otherwise},\\
\end{array} \right.
$$
which leads to (3).

Similarly, (4) holds. $\qed$

Let $C_{i,j}^l(v)$ be a matrix with rows indexed by ${V \choose i}$
and columns indexed by $V \choose j$, whose
 $yz$-entry is defined by
$$
(C_{i,j}^l(v))_{yz}={|y\cap z| \choose l}.
$$
Let $H_{i,j}^l(v)$ be a $(0,1)$-matrix whose rows and columns are
indexed by elements of ${V \choose i}$ and ${V \choose j}$,
respectively, and the $yz$-entry is defined by
$$
(H_{i,j}^l(v))_{yz}=\left\{ \begin{array}{ll}
1,& \textrm{if $|y\cap z|=l$},\\
0,& \textrm{otherwise}.\\
\end{array} \right.
$$
These two matrices may be considered as \emph{intersection matrices}
in the sense that the $yz$-entry only depends on $|y\cap z|$.
Observe $C_{i,j}^0(v)$ is the all-one matrix and
$C_{i,j}^{\min(i,j)}(v)=W_{i,j}(v)$ $(i\leq j)$ or
$(W_{j,i}(v))^\textrm{t}$ $(i>j)$. We adopt the convention that
$C_{i,j}^l(v)=0$ for any integer $l$ such that $l<0$ or
$l>\min(i,j)$. Note that
\begin{eqnarray}\label{17}
C_{i,j}^l(v)=\sum_{g=l}^{\min(i,j)}{g \choose l}H_{i,j}^g(v).
\end{eqnarray}

\begin{lemma}\label{8}
Let $V$ be a set of size $v$. Write $W_{i,j}=W_{i,j}(v)$ and
$C_{i,j}^l=C_{i,j}^l(v)$. Then

\emph{(i)} $W_{i,j}^\emph{t}W_{i,k}=C_{j,k}^i$.

\emph{(ii)} $C_{i,j}^lW_{j,k}={k-l \choose j-l}C_{i,k}^l$.

\emph{(iii)}
$W_{i,k}W_{j,k}^\emph{t}=\sum\limits_{l=\max(0,i+j-k)}^{\min(i,j)}{v-i-j
\choose k-i-j+l}C_{i,j}^l$.

\emph{(iv)}
$W_{i,j}C_{j,k}^l=\sum\limits_{h=\max(0,l+j-i)}^{\min(l,i)}{v-l-i
\choose j-l-i+h}{k-h \choose l-h}C_{i,k}^h$.

\emph{(v)}
$C_{i,j}^lC_{j,k}^s=\sum\limits_{h=\max(0,l+s-j)}^{\min(l,s)}{v-l-s
\choose j-l-s+h}{i-h \choose l-h}{k-h \choose s-h}C_{i,k}^h$.
\end{lemma}
\textbf{Proof.} (i) See \cite  {mm}.

(ii) Immediate from (\ref{3}) and (i).

(iii) We claim that
\begin{eqnarray}\label{20}
W_{i,i+1}W_{j,i+1}^\textrm{t}=(v-i-j)C_{i,j}^j+C_{i,j}^{j-1} \
(j\leq i+1).
\end{eqnarray}
When $j=i+1$, $W_{i,i+1}W_{j,i+1}^\textrm{t}=W_{i,j}=C_{i,j}^{j-1}$,
(\ref{20}) holds. We now assume $j\leq i$. For any $y\in{V \choose
i}$ and $z\in{V \choose j}$,
\begin{eqnarray*}
&&(W_{i,i+1}W_{j,i+1}^\textrm{t})_{yz}\\
&=&\sum_{w\in{V \choose i+1}}(W_{i,i+1})_{yw}(W_{j,i+1}^\textrm{t})_{wz}\\
&=&|\{w\mid (y\cup z)\subseteq w,w\in{V \choose i+1}\}|\\
&=&\left\{ \begin{array}{lll} v-i,& |y\cap z|=j,\\
1,& |y\cap z|=j-1,\\
0,& |y\cap z|\leq j-2,\\
\end{array} \right.
\end{eqnarray*}
which implies (\ref{20}).

Next we show that
\begin{eqnarray}\label{34}
W_{i,k}W_{j,k}^\textrm{t}=\sum_{s=0}^{\min(k-i,j)}{v-i-j \choose
k-i-s}C_{i,j}^{j-s}.
\end{eqnarray}
Observe (\ref{34}) holds when $i=k$. By induction, (\ref{3}),
(\ref{20}) and (i),
\begin{eqnarray*}
&&W_{i-1,k}W_{j,k}^\textrm{t}\\
&=&\frac{1}{k-i+1}W_{i-1,i}W_{i,k}W_{j,k}^\textrm{t}\\
&=&\frac{1}{k-i+1}\sum_{s=0}^{\min(k-i,j)}{v-i-j \choose
k-i-s}W_{i-1,i}C_{i,j}^{j-s}\\
&=&\sum_{s=0}^{\min(k-i,j)}{v-i-j
\choose
k-i-s}\left(\frac{v-i-j+s+1}{k-i+1}C_{i-1,j}^{j-s}+\frac{s+1}{k-i+1}C_{i-1,j}^{j-s-1}\right)\\
&=&\sum_{s=0}^{\min(k-i+1,j)}{v-i-j+1 \choose
k-i-s+1}C_{i-1,j}^{j-s}.
\end{eqnarray*}
Then (\ref{34}) is obtained by induction, concluding (iii).

(iv) Immediate from (i), (ii) and (iii).

(v) Obtained by (i), (ii) and (iv).   $\qed$

\section{The Terwilliger algebra}

In this section we fix $x\in {\Omega \choose m}$, then consider the
Terwilliger algebra $\mathcal {T}=\mathcal {T}(x)$ of $J(n,m,m+1)$
when $n\geq 3m$.

For $0\leq i,j\leq 2m+1$, any matrix $M$ indexed by elements in
$\Gamma_i(x)\times\Gamma_j(x)$ can be embedded into
Mat$_X(\mathbb{C})$ by
$$
L(M)_{\Gamma_k(x)\times\Gamma_l(x)}=\left\{ \begin{array}{ll}
M,& \textrm{if $k=i$ and $l=j$},\\
0,& \textrm{otherwise}.\\
\end{array} \right.
$$
For $0\leq i,j\leq 2m+1$, let
\begin{eqnarray*}
&&\mathcal {M}_{i,j}\\
&=&\textrm{Span}\{C_{m-\lfloor\frac{i}{2}\rfloor,m-\lfloor\frac{j}{2}\rfloor}^l(m)\otimes
C_{\lceil\frac{i}{2}\rceil,\lceil\frac{j}{2}\rceil}^s(n-m),\\
{}&&0\leq
l\leq\min(m-\lfloor\frac{i}{2}\rfloor,m-\lfloor\frac{j}{2}\rfloor),
\ 0\leq s\leq\min(\lceil\frac{i}{2}\rceil,\lceil\frac{j}{2}\rceil)
\}.
\end{eqnarray*}
Let
\begin{eqnarray}\label{9}
\mathcal {M}=\bigoplus\limits_{i,j=0}^{2m+1}L(\mathcal {M}_{i,j}),
\end{eqnarray}
where $L(\mathcal {M}_{i,j})=\{L(M)\mid M\in \mathcal {M}_{i,j}\}$.

Note that $\mathcal {M}$ is a vector space. By Lemma \ref{8} (v) we
have $\mathcal {M}$ is an algebra. In the remaining of this section
we shall prove $\mathcal {T}=\mathcal {M}$.

We begin with a lemma.

\begin{lemma}\label{11}
The Terwilliger algebra $\mathcal {T}$ is a subalgebra of $\mathcal
{M}$.
\end{lemma}
\textbf{Proof.} By Lemma \ref{4} we have $ A\in \mathcal {M}$. For
$0\leq i\leq 2m+1$, since
$$
E_i^*=E_i^*(x)=L(C_{m-\lfloor\frac{i}{2}\rfloor,m-\lfloor\frac{i}{2}\rfloor}^{m-\lfloor\frac{i}{2}\rfloor}(m)\otimes
C_{\lceil\frac{i}{2}\rceil,\lceil\frac{i}{2}\rceil}^{\lceil\frac{i}{2}\rceil}(n-m))\in
\mathcal {M},
$$
we get $\mathcal {T}\subseteq\mathcal {M}$. $\qed$

For $0\leq i,j\leq 2m+1$, let $\mathcal {T}_{i,j}=\{M_{i,j}\mid
M\in\mathcal {T}\}$, where $M_{i,j}$ is the submatrix of $M$ with
rows indexed by $\Gamma_i(x)$ and columns indexed by $\Gamma_j(x)$.
Since $\mathcal {T}$ is an algebra, each $\mathcal {T}_{i,j}$ is a
linear space. Since $\mathcal {T}E_j^*\mathcal {T}\subseteq\mathcal
{T}$, $(\mathcal {T}E_j^*\mathcal {T})_{i,k}\subseteq\mathcal
{T}_{i,k}$, which gives
\begin{eqnarray}\label{23}
\mathcal {T}_{i,j}\mathcal {T}_{j,k}\subseteq\mathcal {T}_{i,k}.
\end{eqnarray}
Since $A$, $E_i^*\in\mathcal {T}$, we have
$AE_{i_2}^*AE_{i_3}^*\cdots AE_{i_{p-1}}^*A\in\mathcal {T}$, which
follows that
\begin{eqnarray}\label{13}
A_{i_1,i_2}A_{i_2,i_3}\cdots A_{i_{p-2},i_{p-1}}A_{i_{p-1},i_p}\in
\mathcal {T}_{i_1,i_p},
\end{eqnarray}
where $0\leq i_1,i_2,\ldots,i_p\leq 2m+1$.

\begin{lemma}\label{12}
For $2i+2\leq j\leq 2m+1$ and $0\leq s\leq i+1$, we have
$$
C_{m-i-1,m-\lfloor\frac{j}{2}\rfloor}^{m-\lfloor\frac{j}{2}\rfloor}(m)\otimes
C_{i+1,\lceil\frac{j}{2}\rceil}^s(n-m)\in\mathcal {T}_{2i+2,j}.
$$
\end{lemma}
\textbf{Proof.} We use induction on $s$ ($s$ decreasing from $i+1$
to $0$).

By (\ref{13}), for $j> 2i+2$ we have
$A_{2i+2,2i+3}A_{2i+3,2i+4}\cdots A_{j-1,j}\in\mathcal
{T}_{2i+2,j}$, which yields that
\begin{eqnarray}\label{15}
C_{m-i-1,m-\lfloor\frac{j}{2}\rfloor}^{m-\lfloor\frac{j}{2}\rfloor}(m)\otimes
C_{i+1,\lceil\frac{j}{2}\rceil}^{i+1}(n-m)\in\mathcal {T}_{2i+2,j}.
\end{eqnarray}
When $j=2i+2$ we pick $I_{m \choose m-i-1}\otimes I_{n-m \choose
i+1}\in\mathcal {T}_{2i+2,2i+2}$, which also satisfies (\ref{15}).

Assume that
$C_{m-i-1,m-\lfloor\frac{j}{2}\rfloor}^{m-\lfloor\frac{j}{2}\rfloor}(m)\otimes
C_{i+1,\lceil\frac{j}{2}\rceil}^s(n-m)\in\mathcal {T}_{2i+2,j}$. By
(\ref{23}) and (\ref{13}) we obtain
\begin{eqnarray}\label{35}
(C_{m-i-1,m-\lfloor\frac{j}{2}\rfloor}^{m-\lfloor\frac{j}{2}\rfloor}(m)\otimes
C_{i+1,\lceil\frac{j}{2}\rceil}^s(n-m))(A_{j,j+1}A_{j+1,j})\in
\mathcal {T}_{2j+2,j}\mathcal {T}_{j,j}\subseteq\mathcal
{T}_{2i+2,j},
\end{eqnarray}
\begin{eqnarray}\label{33}
(C_{m-i-1,m-\lfloor\frac{j}{2}\rfloor}^{m-\lfloor\frac{j}{2}\rfloor}(m)\otimes
C_{i+1,\lceil\frac{j}{2}\rceil}^s(n-m))(A_{j,j-1}A_{j-1,j})\in
\mathcal {T}_{2j+2,j}\mathcal {T}_{j,j}\subseteq\mathcal
{T}_{2i+2,j}.
\end{eqnarray}
When $j$ is even, by Lemma \ref{4}, Lemma \ref{8} (iii) and (v),
(\ref{35}) leads to
$$
aC_{m-i-1,m-\frac{j}{2}}^{m-\frac{j}{2}}(m)\otimes
C_{i+1,\frac{j}{2}}^s(n-m)+bC_{m-i-1,m-\frac{j}{2}}^{m-\frac{j}{2}}(m)\otimes
C_{i+1,\frac{j}{2}}^{s-1}(n-m)\in\mathcal {T}_{2i+2,j},
$$
where $a=(n-m-s-\frac{j}{2})(\frac{j}{2}-s+1)$ and
$b=(i-s+2)(\frac{j}{2}-s+1)$. Similarly when $j$ is odd, (\ref{33})
yields that
$$
a'C_{m-i-1,m-\lfloor\frac{j}{2}\rfloor}^{m-\lfloor\frac{j}{2}\rfloor}(m)\otimes
C_{i+1,\lceil\frac{j}{2}\rceil}^s(n-m)+b'C_{m-i-1,m-\lfloor\frac{j}{2}\rfloor}^{m-\lfloor\frac{j}{2}\rfloor}(m)\otimes
C_{i+1,\lceil\frac{j}{2}\rceil}^{s-1}(n-m)\in\mathcal {T}_{2i+2,j},
$$
where
$a'=(n-m-s-\lceil\frac{j}{2}\rceil+1)(\lceil\frac{j}{2}\rceil-s)$
and $b'=(i-s+2)(\lceil\frac{j}{2}\rceil-s+1)$. Since $s\leq i+1\leq
\lceil\frac{j}{2}\rceil$, $(i-s+2)(\lceil\frac{j}{2}\rceil-s+1)\neq
0$. Thus we have
$C_{m-i-1,m-\lfloor\frac{j}{2}\rfloor}^{m-\lfloor\frac{j}{2}\rfloor}(m)\otimes
C_{i+1,\lceil\frac{j}{2}\rceil}^{s-1}(n-m)\in\mathcal {T}_{2i+2,j}$.

Hence the desired result follows. $\qed$

\begin{lemma}\label{14}
The algebra $\mathcal {M}$ is a subalgebra of $\mathcal {T}$.
\end{lemma}
\textbf{Proof.} During this proof we will omit the symbol $(m)$ from
matrices in front of ``$\otimes$'', and omit $(n-m)$ from matrices
behind ``$\otimes$''.

In order to get the desired conclusion, we only need to show that
$\mathcal {M}_{i,j}\subseteq\mathcal {T}_{i,j}$ for $0\leq i,j\leq
2m+1$. Write $\mathcal {M}_{i,j}^\textrm{t}=\{M^\textrm{t}\mid
M\in\mathcal {M}_{i,j}\}$ and $\mathcal
{T}_{i,j}^\textrm{t}=\{M^\textrm{t}\mid M\in\mathcal {T}_{i,j}\}$  .
Since $\mathcal {M}_{j,i}=\mathcal {M}_{i,j}^\textrm{t}$ and
$\mathcal {T}_{j,i}=\mathcal {T}_{i,j}^\textrm{t}$, it suffices to
prove $\mathcal {M}_{i,j}\subseteq\mathcal {T}_{i,j}$ for $i\leq j$.
We use induction on $i$.

\medskip
\textbf{Step 1.} Show $\mathcal {M}_{0,j}\subseteq\mathcal
{T}_{0,j}$ $(0\leq j\leq 2m+1)$.

According to (\ref{9}),
$$
\mathcal
{M}_{0,j}=\textrm{Span}\{C_{m,m-\lfloor\frac{j}{2}\rfloor}^l\otimes
C_{0,\lceil\frac{j}{2}\rceil}^0, \ \ 0\leq l\leq
m-\lfloor\frac{j}{2}\rfloor\}.
$$
For any $l\in\{0,1,\ldots,m-\lfloor\frac{j}{2}\rfloor\}$,
$$
C_{m,m-\lfloor\frac{j}{2}\rfloor}^l\otimes
C_{0,\lceil\frac{j}{2}\rceil}^0={m-\lfloor\frac{j}{2}\rfloor \choose
l}W_{m-\lfloor\frac{j}{2}\rfloor,m}^\textrm{t}\otimes
W_{0,\lceil\frac{j}{2}\rceil},
$$
while by (\ref{3}) and Lemma \ref{4},
$$
A_{0,1}A_{1,2}\cdots
A_{j-1,j}=\lfloor\frac{j}{2}\rfloor!\lceil\frac{j}{2}\rceil!W_{m-\lfloor\frac{j}{2}\rfloor,m}^\textrm{t}\otimes
W_{0,\lceil\frac{j}{2}\rceil}.
$$
Hence we get $\mathcal {M}_{0,j}\subseteq\mathcal {T}_{0,j}$ from
(\ref{13}).

\medskip
\textbf{Step 2.} Assume that $\mathcal {M}_{p,j}\subseteq\mathcal
{T}_{p,j}$ for $p\leq 2i$. We will show that $\mathcal
{M}_{2i+1,j}\subseteq\mathcal {T}_{2i+1,j}$ and $\mathcal
{M}_{2i+2,j}\subseteq\mathcal {T}_{2i+2,j}$.

\medskip
\textbf{Step 2.1.} Show $\mathcal {M}_{2i+1,j}\subseteq\mathcal
{T}_{2i+1,j}$ $(2i+1\leq j\leq 2m+1)$.

It suffices to prove
\begin{eqnarray}\label{31}
C_{m-i,m-\lfloor\frac{j}{2}\rfloor}^l\otimes
C_{i+1,\lceil\frac{j}{2}\rceil}^s\in\mathcal {T}_{2i+1,j},
\end{eqnarray}
where $0\leq l\leq m-\lfloor\frac{j}{2}\rfloor$, $0\leq s\leq i+1$.

By inductive hypothesis,
$$
C_{m-i,m-\lfloor\frac{j}{2}\rfloor}^l\otimes
C_{i,\lceil\frac{j}{2}\rceil}^s\in\mathcal
{M}_{2i,j}\subseteq\mathcal {T}_{2i,j}, \quad 0\leq l\leq
m-\lfloor\frac{j}{2}\rfloor,\ 0\leq s\leq i.
$$
Since
$$
A_{2i,2i+1}^\textrm{t}=I_{m \choose m-i}\otimes
W_{i,i+1}^\textrm{t}\in\mathcal
{M}_{2i,2i+1}^\textrm{t}\subseteq\mathcal {T}_{2i,2i+1}^\textrm{t},
$$
we have
$$
(I_{m \choose m-i}\otimes
W_{i,i+1}^\textrm{t})(C_{m-i,m-\lfloor\frac{j}{2}\rfloor}^l\otimes
C_{i,\lceil\frac{j}{2}\rceil}^s)\in\mathcal
{T}_{2i,2i+1}^\textrm{t}\mathcal {T}_{2i,j}\subseteq\mathcal
{T}_{2i+1,j},
$$
which by Lemma \ref{8} (ii) follows that (\ref{31}) holds for $0\leq
l\leq m-\lfloor\frac{j}{2}\rfloor$, $0\leq s\leq i$.

By (\ref{13}), for $j>2i+1$ we get
$$
A_{2i+1,2i+2}A_{2i+2,2i+3}\cdots A_{j-1,j}\in\mathcal {T}_{2i+1,j},
$$
which yields that
\begin{eqnarray}\label{24}
W_{m-\lfloor\frac{j}{2}\rfloor,m-i}^\textrm{t}\otimes
W_{i+1,\lceil\frac{j}{2}\rceil}=C_{m-i,m-\lfloor\frac{j}{2}\rfloor}^{m-\lfloor\frac{j}{2}\rfloor}\otimes
C_{i+1,\lceil\frac{j}{2}\rceil}^{i+1}\in\mathcal {T}_{2i+1,j}.
\end{eqnarray}
When $j=2i+1$ we pick $I_{m \choose m-i}\otimes I_{n-m \choose
i+1}\in\mathcal {T}_{2i+1,2i+1}$, which also satisfies (\ref{24}).

\emph{Case} 1. $j=2m+1$ or $2m$.

In this case, (\ref{24}) implies that $C_{m-i,0}^0\otimes
C_{i+1,\lceil\frac{j}{2}\rceil}^{i+1}\in\mathcal {T}_{2i+1,j}$,
which means (\ref{31}) holds for $l=0$ and $s=i+1$.

\emph{Case} 2. $j\leq 2m-1$.

For $j+1\leq k\leq 2m$, let
$$
N_{j,k}=A_{j,j+1}A_{j+1,j+2}\cdots A_{k-1,k}.
$$
Again by (\ref{13})
$$
A_{2i+1,2i+2}A_{2i+2,2i+3}\cdots
A_{j-1,j}N_{j,k}N_{j,k}^\textrm{t}\in\mathcal {T}_{2i+1,j}.
$$
By (\ref{3}), Lemma \ref{4} and Lemma \ref{8}(i), (iii), we obtain
\begin{eqnarray}\label{16}
cC_{m-i,m-\lfloor\frac{j}{2}\rfloor}^{m-\lfloor\frac{k}{2}\rfloor}\otimes
\left(\sum_{h=\max(0,i+1+\lceil\frac{j}{2}\rceil-\lceil\frac{k}{2}\rceil)}^{i+1}{n-m-i-1-\lceil\frac{j}{2}\rceil
\choose
\lceil\frac{k}{2}\rceil-i-1-\lceil\frac{j}{2}\rceil+h}C_{i+1,\lceil\frac{j}{2}\rceil}^h\right)\in\mathcal
{T}_{2i+1,j},
\end{eqnarray}
where
$c=(\lfloor\frac{k}{2}\rfloor-i)!(\lceil\frac{k}{2}\rceil-i-1)!(\lfloor\frac{k}{2}\rfloor-\lfloor\frac{j}{2}\rfloor)!
(\lceil\frac{k}{2}\rceil-\lceil\frac{j}{2}\rceil)!\neq 0$. We
consider the coefficient of
$C_{m-i,m-\lfloor\frac{j}{2}\rfloor}^{m-\lfloor\frac{k}{2}\rfloor}\otimes
C_{i+1,\lceil\frac{j}{2}\rceil}^{i+1}$, which is
$c{n-m-i-1-\lceil\frac{j}{2}\rceil \choose
\lceil\frac{k}{2}\rceil-\lceil\frac{j}{2}\rceil}$. Since $0\leq
2i+1\leq j\leq k-1\leq 2m-1$ and $n\geq 3m$, we get
$$
n-m-i-1-\lceil\frac{j}{2}\rceil\geq
n-m-m-\lceil\frac{j}{2}\rceil\geq
m-\lceil\frac{j}{2}\rceil\geq\lceil\frac{k}{2}\rceil-\lceil\frac{j}{2}\rceil\geq
0,
$$
and so $c{n-m-i-1-\lceil\frac{j}{2}\rceil \choose
\lceil\frac{k}{2}\rceil-\lceil\frac{j}{2}\rceil}\neq 0$. Since
(\ref{31}) holds for $s\in\{0,1,\ldots,i\}$, (\ref{31}) also holds
for $s=i+1$ by (\ref{24}) and (\ref{16}).

\medskip
\textbf{Step 2.2.} Show $\mathcal {M}_{2i+2,j}\subseteq \mathcal
{T}_{2i+2,j}$ $(2i+2\leq j\leq 2m+1)$.

It suffices to prove
\begin{eqnarray}\label{25}
C_{m-i-1,m-\lfloor\frac{j}{2}\rfloor}^l\otimes
C_{i+1,\lceil\frac{j}{2}\rceil}^s\in \mathcal {T}_{2i+2,j}, \quad
0\leq l\leq m-\lfloor\frac{j}{2}\rfloor, \ 0\leq s\leq i+1.
\end{eqnarray}
By the inductive assumption, for $0\leq l\leq
m-\lfloor\frac{j}{2}\rfloor$ and $0\leq s\leq i+1$,
$$
C_{m-i,m-\lfloor\frac{j}{2}\rfloor}^l\otimes
C_{i+1,\lceil\frac{j}{2}\rceil}^s\in\mathcal
{M}_{2i+1,j}\subseteq\mathcal {T}_{2i+1,j}.
$$
Since
$$
A_{2i+1,2i+2}^\textrm{t}=W_{m-i-1,m-i}\otimes I_{n-m \choose
i+1}\in\mathcal {T}_{2i+1,2i+2}^\textrm{t},
$$
by (\ref{23}) we have
\begin{eqnarray}
&&(W_{m-i-1,m-i}\otimes I_{n-m \choose
i+1})(C_{m-i,m-\lfloor\frac{j}{2}\rfloor}^l\otimes
C_{i+1,\lceil\frac{j}{2}\rceil}^s)\nonumber\\
&=&(W_{m-i-1,m-i}C_{m-i,m-\lfloor\frac{j}{2}\rfloor}^l)\otimes
C_{i+1,\lceil\frac{j}{2}\rceil}^s\nonumber\\
&\in&\mathcal {T}_{2i+1,2i+2}^\textrm{t}\mathcal
{T}_{2i+1,j}\nonumber\\
&\subseteq&\mathcal {T}_{2i+2,j}.\label{26}
\end{eqnarray}
By Lemma \ref{8} (iv),
$$
W_{m-i-1,m-i}C_{m-i,m-\lfloor\frac{j}{2}\rfloor}^l=
(i+1-l)C_{m-i-1,m-\lfloor\frac{j}{2}\rfloor}^l+(m-\lfloor\frac{j}{2}\rfloor-l+1)C_{m-i-1,m-\lfloor\frac{j}{2}\rfloor}^{l-1}.
$$
Thus (\ref{26}) leads to
\begin{eqnarray}\label{27}
[(i+1-l)C_{m-i-1,m-\lfloor\frac{j}{2}\rfloor}^l+(m-\lfloor\frac{j}{2}\rfloor-l+1)C_{m-i-1,m-\lfloor\frac{j}{2}\rfloor}^{l-1}]\otimes
C_{i+1,\lceil\frac{j}{2}\rceil}^s\in\mathcal {T}_{2i+2,j},
\end{eqnarray}
where $0\leq l\leq m-\lfloor\frac{j}{2}\rfloor$, $0\leq s\leq i+1$.\
Since the coefficient of
$C_{m-i-1,m-\lfloor\frac{j}{2}\rfloor}^{l-1}\otimes
C_{i+1,\lceil\frac{j}{2}\rceil}^s$ in (\ref{27}) is
$m-\lfloor\frac{j}{2}\rfloor-l+1\neq 0$, by Lemma \ref{12} we get
(\ref{25}).

Hence the desired result follows. $\qed$

\begin{thm}\label{21}
Let $J(n,m,m+1)$ be the incidence graph of Johnson geometry with
$n\geq 3m$. Let $\mathcal {T}=\mathcal {T}(x)$ be the Terwilliger
algebra of $J(n,m,m+1)$ with respect to an $m$-subset $x$ and
$\mathcal {M}$ be the corresponding algebra defined in $(\ref{9})$.
Then $\mathcal {T}=\mathcal {M}$.
\end{thm}
\textbf{Proof.} Combining Lemma \ref{11} and Lemma \ref{14}, the
proof of Theorem \ref{21} is completed. $\qed$

\begin{cor}\label{32}
With reference to Theorem \ref{21} $J(n,m,m+1)$ is thin with respect
to $x$.
\end{cor}
\textbf{Proof.} By Theorem \ref{21} we get
$$
E_i^*\mathcal
{T}E_i^*=\textrm{Span}\{L(C_{m-\lfloor\frac{i}{2}\rfloor,m-\lfloor\frac{i}{2}\rfloor}^l(m)\otimes
C_{\lceil\frac{i}{2}\rceil,\lceil\frac{i}{2}\rceil}^s(n-m)), \ 0\leq
l\leq m-\lfloor\frac{i}{2}\rfloor, \ 0\leq s\leq
\lceil\frac{i}{2}\rceil\},
$$
where $i=0,1,\ldots,D(x)$. Since each element of $E_i^*\mathcal
{T}E_i^*$ is symmetric, we get the conclusion from \cite [Theorem
13]{ter4}. $\qed$

\section{The basis of the Terwilliger algebra}

In this section we shall determine the basis and the dimension of
$\mathcal {T}$.

\begin{thm}\label{18}
Let $G_{i,j}=\{g\mid
H_{m-\lfloor\frac{i}{2}\rfloor,m-\lfloor\frac{j}{2}\rfloor}^g(m)\neq
0\}$, $R_{i,j}=\{r\mid
H_{\lceil\frac{i}{2}\rceil,\lceil\frac{j}{2}\rceil}^r(n-m)\neq 0\}$,
and $\mathcal {T}$ be as in Theorem \ref{21}. Then we have
\begin{eqnarray}\label{2}
\{L(H_{m-\lfloor\frac{i}{2}\rfloor,m-\lfloor\frac{j}{2}\rfloor}^g(m)\otimes
H_{\lceil\frac{i}{2}\rceil,\lceil\frac{j}{2}\rceil}^r(n-m)), \ g\in
G_{i,j}, \ r\in R_{i,j}\}_{i,j=0}^{2m+1}
\end{eqnarray}
as well as
\begin{eqnarray}\label{10}
\{L(C_{m-\lfloor\frac{i}{2}\rfloor,m-\lfloor\frac{j}{2}\rfloor}^l(m)\otimes
C_{\lceil\frac{i}{2}\rceil,\lceil\frac{j}{2}\rceil}^s(n-m)), \ l\in
G_{i,j}, \ s\in R_{i,j}\}_{i,j=0}^{2m+1}
\end{eqnarray}
are two bases of $\mathcal {T}$.
\end{thm}
\textbf{Proof.} Without loss of generality, we assume $i\leq j$.
Since $H_{i,j}^l(v)\neq 0$ if and only if $\max(0,i+j-v)\leq l\leq
\min(i,j)$, we have $\lceil\frac{i}{2}\rceil-|R_{i,j}|+1\leq r \leq
\lceil\frac{i}{2}\rceil$ when $r\in R_{i,j}$. By (\ref{17}) we
obtain
\begin{eqnarray}\label{19}
C_{\lceil\frac{i}{2}\rceil,\lceil\frac{j}{2}\rceil}^r(n-m)
=\sum_{h=r}^{\lceil\frac{i}{2}\rceil}{h \choose
r}H_{\lceil\frac{i}{2}\rceil,\lceil\frac{j}{2}\rceil}^h(n-m),
\end{eqnarray}
which implies that
$H_{\lceil\frac{i}{2}\rceil,\lceil\frac{j}{2}\rceil}^r(n-m)$ $(r\in
R_{i,j})$ is a linear combination of
$\{C_{\lceil\frac{i}{2}\rceil,\lceil\frac{j}{2}\rceil}^s(n-m)\}_{s\in
R_{i,j}}$. Similarly,
$H_{m-\lfloor\frac{i}{2}\rfloor,m-\lfloor\frac{j}{2}\rfloor}^g(m)$
$(g\in G_{i,j})$ can be expressed as a linear combination of
$\{C_{m-\lfloor\frac{i}{2}\rfloor,m-\lfloor\frac{j}{2}\rfloor}^l(m)\}_{l\in
G_{i,j}}$. Hence
$$
H_{m-\lfloor\frac{i}{2}\rfloor,m-\lfloor\frac{j}{2}\rfloor}^g(m)\otimes
H_{\lceil\frac{i}{2}\rceil,\lceil\frac{j}{2}\rceil}^r(n-m)\in
\mathcal {M}_{i,j}.
$$
Again by (\ref{17}), for $0\leq l\leq m-\lfloor\frac{j}{2}\rfloor$
and $0\leq s\leq \lceil\frac{i}{2}\rceil$,
\begin{eqnarray*}
&&C_{m-\lfloor\frac{i}{2}\rfloor,m-\lfloor\frac{j}{2}\rfloor}^l(m)\otimes
C_{\lceil\frac{i}{2}\rceil,\lceil\frac{j}{2}\rceil}^s(n-m)\\
&=&\sum_{g=l}^{m-\lfloor\frac{j}{2}\rfloor}{g \choose
l}H_{m-\lfloor\frac{i}{2}\rfloor,m-\lfloor\frac{j}{2}\rfloor}^g(m)\otimes
\sum_{r=s}^{\lceil\frac{i}{2}\rceil}{r \choose
s}H_{\lceil\frac{i}{2}\rceil,\lceil\frac{j}{2}\rceil}^r(n-m).
\end{eqnarray*}
Observe that
$H_{m-\lfloor\frac{i}{2}\rfloor,m-\lfloor\frac{j}{2}\rfloor}^g(m)\otimes
H_{\lceil\frac{i}{2}\rceil,\lceil\frac{j}{2}\rceil}^r(n-m)$ $(g\in
G_{i,j}, r\in R_{i,j})$ are linearly independent. Then
$\{H_{m-\lfloor\frac{i}{2}\rfloor,m-\lfloor\frac{j}{2}\rfloor}^g(m)\otimes
H_{\lceil\frac{i}{2}\rceil,\lceil\frac{j}{2}\rceil}^r(n-m)\}_{g\in
G_{i,j}, r\in R_{i,j}}$ is a basis of $\mathcal {M}_{i,j}$.
Therefore (\ref{2}) is a basis of $\mathcal {T}$.

Furthermore by (\ref{19}) we can get
$\{C_{m-\lfloor\frac{i}{2}\rfloor,m-\lfloor\frac{j}{2}\rfloor}^l(m)\otimes
C_{\lceil\frac{i}{2}\rceil,\lceil\frac{j}{2}\rceil}^s(n-m)\}_{l\in
G_{i,j}, s\in R_{i,j}}$ is also a basis of $\mathcal {M}_{i,j}$,
which follows that (\ref{10}) is a basis of $\mathcal {T}$.

This ends our proof. $\qed$

\begin{cor}\label{28}
With reference to Theorem \ref{21} we get the dimension of $\mathcal
{T}$ is
$$
\dim \mathcal {T}=\left\{ \begin{array}{lll} \frac{1}{12}(m+1)(m+2)(m+3)(3m+10)-4,& \emph{if $n=3m$},\\
\frac{1}{12}(m+1)(m+2)(m+3)(3m+10)-1,& \emph{if $n=3m+1$},\\
\frac{1}{12}(m+1)(m+2)(m+3)(3m+10),& \emph{if $n\geq3m+2$}.\\
\end{array} \right.
$$
\end{cor}
\textbf{Proof.} By Theorem \ref{18} we get
\begin{eqnarray*}
\dim \mathcal {T}&=&\sum_{i,j=0}^{2m+1}|G_{i,j}||R_{i,j}|\\
&=&\sum_{i,j=0}^{2m+1}(\min(m-\lfloor\frac{i}{2}\rfloor,m-\lfloor\frac{j}{2}\rfloor)
-\max(0,m-\lfloor\frac{i}{2}\rfloor-\lfloor\frac{j}{2}\rfloor)+1)\\
&&\cdot(\min(\lceil\frac{i}{2}\rceil,\lceil\frac{j}{2}\rceil)
-\max(0,\lceil\frac{i}{2}\rceil+\lceil\frac{j}{2}\rceil-n+m)+1).
\end{eqnarray*}
By zigzag calculation, we get the desired result.  $\qed$

\section{Concluding Remark}

We conclude this paper with the following remarks:

(i) Let $J(n,m)$ be the Johnson graph with $n\geq 3m$. Fix a vertex
$x$ of $J(n,m)$. Let $\mathcal {T}'=\mathcal {T}'(x)$ and $\mathcal
{T}=\mathcal {T}(x)$ be the Terwilliger algebra of $J(n,m)$ and
$J(n,m,m+1)$ with respect to $x$, respectively. Since
$\bigoplus_{i,j=0}^mE_{2i}^*(x)\mathcal {T}E_{2j}^*(x)$ is an
algebra, $\{L(H_{m-i,m-j}^g(m)\otimes H_{i,j}^r(n-m)),\ g\in
G_{2i,2j}, \ r\in R_{2i,2j}\}_{i,j=0}^{m}$ is a basis of
$\bigoplus_{i,j=0}^mE_{2i}^*(x)\mathcal {T}E_{2j}^*(x)$ by Theorem
\ref{18}. By \cite [Definition 4.2, Lemma 4.4, Theorem 5.9]{johnson}
this basis coincides with that of $\mathcal {T}'$, which implies
that $\mathcal {T}'\simeq\bigoplus_{i,j=0}^mE_{2i}^*(x)\mathcal
{T}E_{2j}^*(x)$.

(ii) Using the same method, the Terwilliger algebra of $J(n,m,m+1)$
with respect to an $(m+1)$-subset may be determined.

\section*{Acknowledgement} The authors would like to thank Professor
Hiroshi Suzuki for proposing this problem and for his many helpful
suggestions. This research is partially supported by    NSF of China
(10871027), NCET-08-0052,   and   the Fundamental Research Funds for
the Central Universities of China.

\end{CJK*}

\end{document}